\documentclass{article}
 \usepackage[T1,T2A]{fontenc}
\usepackage[unicode,hidelinks]{hyperref}
\usepackage{amsmath,amssymb,amsthm,orcidlink,booktabs,pdflscape}
\usepackage[margin=1in]{geometry}
\usepackage{CJKutf8}
\usepackage[final]{microtype}
\title{Cosmological polytopes are dual to order polytopes}

\author{Stefan Forcey\textsuperscript{\orcidlink{0000-0003-2251-5710}}\footnote{Department of Mathematics, University of Akron, Akron, Ohio , USA}\and Ross Glew\textsuperscript{\orcidlink{0000-0002-8100-928X}}\footnote{School of Natural Sciences, Institute for Advanced Study, Princeton, New Jersey, USA}\footnotemark[1]
\\[1em]
\href{mailto:sforcey@uakron.edu}{\texttt{sforcey@uakron.edu}}
\hskip15em
\href{mailto:r.glew@herts.ac.uk}{\texttt{rglew@ias.edu}}
}

\usepackage[french,german,portuguese,italian,latin,russian,ukrainian,british]{babel}
\usepackage[osf]{newpxtext}
\usepackage[varbb,slantedGreek]{newpxmath}
\usepackage[artemisia]{textgreek}
\usepackage[utf8]{inputenc}
\usepackage{tikz-cd}
\usetikzlibrary{babel}
\newcommand{\precdot}{\prec\mathrel{\mkern-5mu}\mathrel{\cdot}}

\begin{document}
\maketitle
\newtheorem{theorem}{Theorem}
\newtheorem{lemma}{Lemma}
\newtheorem{conjecture}{Conjecture}
\newtheorem{corollary}{Corollary}
\theoremstyle{definition}
\newtheorem{defn}{Definition}
\theoremstyle{remark}
\newtheorem{remark}{Remark}
\newtheorem{example}{Example}
\newcommand\cyrillic[1]{\fontfamily{Domitian-TOsF}\selectfont \foreignlanguage{russian}{#1}}

\baselineskip=17pt

\abstract{ We show that every cosmological polytope is combinatorially equivalent to the polar dual of an order polytope, and that the faces of the scattering facet of that cosmological polytope are partitioned into hypergraphic polytopes. As a first consequence, we see a new labeling of any hypergraphic polytope with faces of the order polytope, which displays that hypergraphic polytope as a cross section of the polar dual of that order polytope. A second corollary is that the cosmological polytope is a root polytope, and we describe for it an explicit root polytope realization. We use that realization to show our cross sectional realization of hypergraphic polytopes, as a specialization of Postnikov's cross sectional realization of generalized permutahedra. Our description  allows for natural generalizations of the hypergraphic and cosmological polytopes to pseudo-hypergraphs.  Our corollaries imply formulas for the $f$-vectors and polynomials of rank-2 order polytopes in terms of hypergraphic polytopes.}

\begin{figure}
    \centering
    \includegraphics[width=\linewidth]{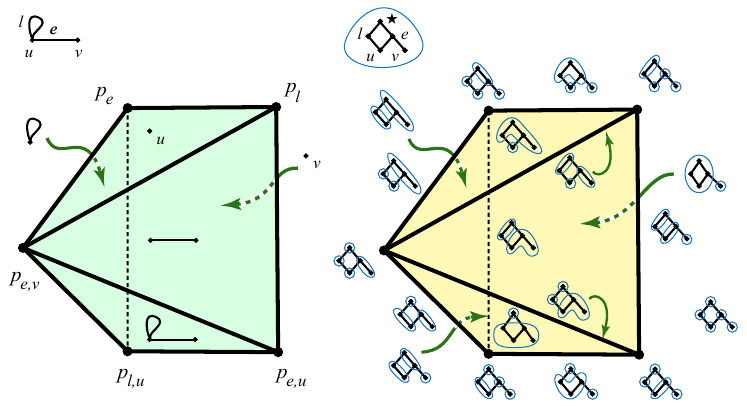}
    \caption{Left: the cosmological polytope   $\mathcal{C}_G$ associated to the pseudograph $G$ from Example~\ref{uno}. Right: the polar dual of the order polytope of the poset $P'(G)$, denoted $\mathcal{O}(P'(G))^{\Delta}$; the face lattice labeled by connected, acyclic partitions of $P'(G)$, as listed in Example~\ref{dos} and demonstrating Theorem~\ref{big}. }
    \label{fig:compare}
\end{figure}

\begin{figure}
     \centering
     \includegraphics[width=\linewidth]{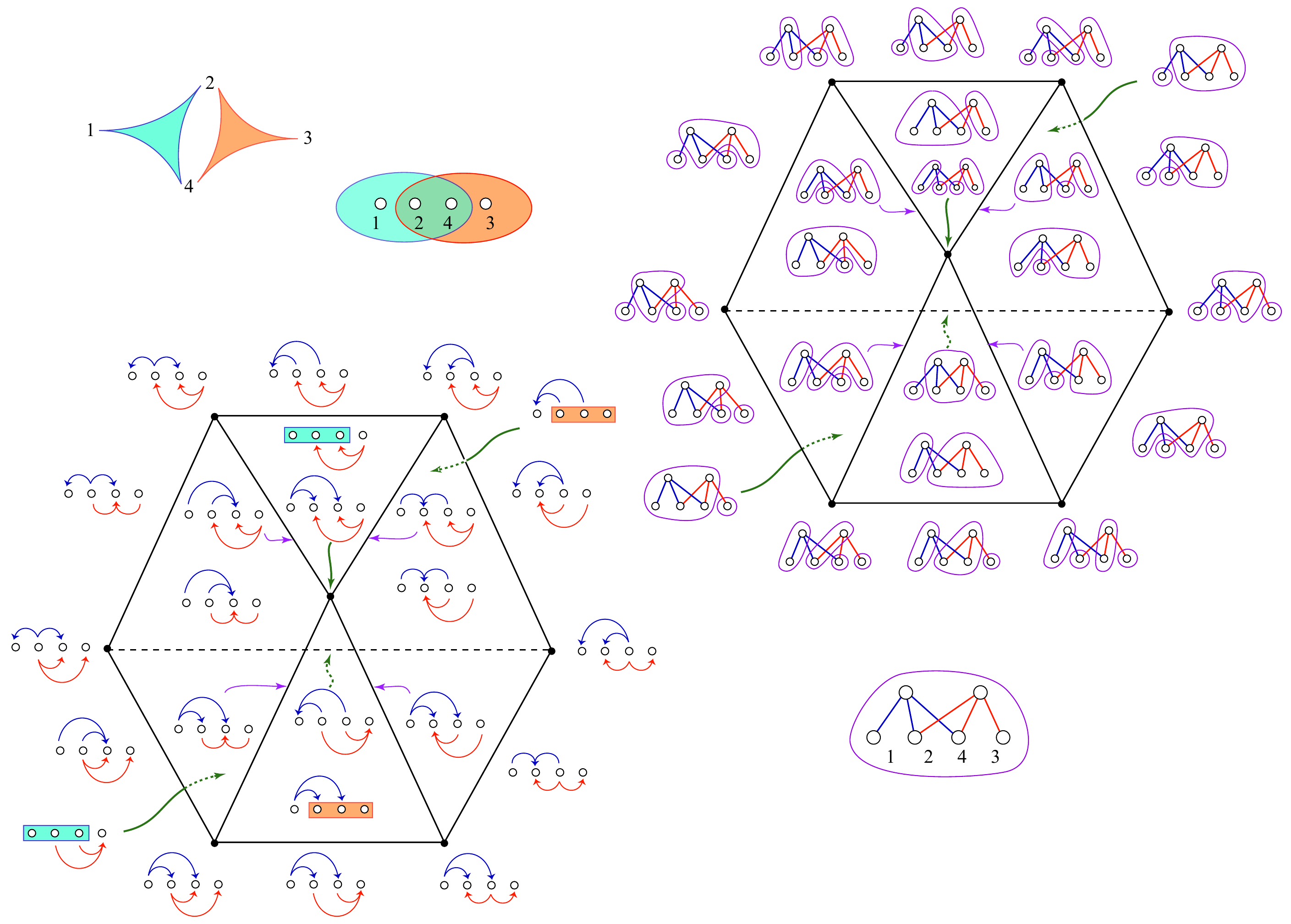}
     \caption{Example of the bijection in the proof of Theorem~\ref{conj}. The hypergraph $G$ pictured two ways at far left has hyperedges $E = \{\{1,2,4\},\{2,3,4\}\}$. The source sets of each hyperedge in the acyclic orientations are the starting nodes for each pair of same color arrows. When the source set is the entire hyperedge the nodes are highlighted. To the right is seen the corresponding labeling with acyclic connected partitions of $P(G)$, which is shown to be realized as a cross section of the polar dual of the order polytope.}
     \label{fig:hyp3ord}
 \end{figure}


\section{Overview}
Cosmological polytopes are a graph-indexed family of polytopes introduced to geometrize individual Feynman graph contributions to the \emph{wavefunction of the universe} \cite{Arkani-Hamed:2}. They have since been studied in their own right for their rich combinatorial structure \cite{juhnke,kuhne,birkemeyer2026cosmological}. Posted concurrently with our present paper is a new study of cosmological polytopes and the emergence of time: \cite{etimes}.

In this paper, we explore connections between cosmological polytopes and three other families: order polytopes, hypergraphic polytopes, and root polytopes. Our first main result, presented in Theorem~\ref{big}, is that every cosmological polytope is combinatorially equivalent to the polar dual of an order polytope associated with a poset constructed from the pseudograph. We establish this equivalence through explicit bijections between vertices and facets that preserve their incidence relations, and hence determine the entire face lattice. This allows combinatorial results for either family of polytopes to be transferred to the other. In particular, combining our result with Theorem D of \cite{williams}, which relates Stanley's order polytope to a root polytope associated with the Hasse diagram of a poset, gives a new realization of the cosmological polytope as a root polytope, presented in Corollary~\ref{newg} of Section~\ref{root}.

Hypergraphic polytopes were recently defined and studied in \cite{berg,pilaud_berg}. Our second main result, presented in Theorem~\ref{conj}, is that every hypergraphic polytope can be realized as a cross-section of the polar dual of an order polytope. We establish this through a bijection between the faces of the hypergraphic polytope and a certain collection of faces of the order polytope. We then use the realization of a facet of the cosmological polytope as a root polytope to give a convenient description of this cross-section.
More generally, we show that the faces of the order polytope of any rank-2 poset can be partitioned into collections whose face posets are equivalent to the dual face posets of (sub-)hypergraphic polytopes. This yields a formula for the $f$-vector of the order polytope in terms of those of hypergraphic polytopes, given in Corollary~\ref{effy}.

The paper is organized as follows. In Sections~\ref{sec:hyp}, \ref{sec:cosmo}, and \ref{sec:order}, we review hypergraphic, cosmological, and order polytopes, respectively. In Section~\ref{sec:equivo}, we establish the combinatorial equivalence between cosmological polytopes and the polar duals of order polytopes associated with graph incidence posets. In Section~\ref{root}, we use the connection to root polytopes to obtain new realizations of cosmological polytopes. Finally, in Section~\ref{sec:hyp2}, we study sections of cosmological polytopes and their relation to hypergraphic polytopes.

\section{Hypergraphic polytopes}\label{sec:hyp}

Throughout we will deal with hypergraphs, usually allowing loops and multiedges as well. To emphasize that, we often refer to pseudo-hypergraphs $G=(V,E)$ where $V$ is the set of nodes, and the collection of hyperedges $E$ is a multiset of subsets of $V$. Simple hypergraphs in contrast have no loops or multiedges.  

Hypergraphic polytopes are defined in \cite{berg} and \cite{pilaud_berg}. Given a connected simple hypergraph $G$, a partial orientation is a selection of a non-empty subset of each hyperedge as the source set, (or head), with the complement in that hyperedge called the target set, (or tail), possibly empty).

\begin{defn}
A (partial) orientation $\mathcal{S}=\{s_e  \}_{e \in E}$ of $G$ is a choice of non-empty source set $\emptyset\neq s_e \subseteq e$ for each hyperedge $e \in E$. We refer to the complements $t_e =e \setminus h_e$ as target sets. 
\end{defn}

\begin{defn}
Given an orientation $\mathcal{S}$ of a hypergraph $G$, we define a directed graph $G/\mathcal{S}$ as follows. Its vertex set $V/\mathcal{S}$ is the set of equivalence classes generated by the relation $s \sim s'$ whenever $s,s'\in s_e$ for some $e\in E$. For each $e\in E$ and each $t\in t_e$, there is a directed edge $[s_e]\longrightarrow [t]$, where $[s_e]$ denotes the common equivalence class containing the elements of $s_e$.
\end{defn}

\begin{defn}
An orientation $\mathcal{S}$ of $G$ is \emph{acyclic} if the graph $G/\mathcal{S}$ contains no directed cycles.
\end{defn}

Acyclic implies that there is no cycle formed by a series of non-empty subsets of sources: $a_1,a_2,\dots,a_k$ with $a_i\subseteq s_i $,  and $a_{i+1} \subseteq t_i$, and $a_k\subseteq s_1.$ 

\begin{figure}
    \centering
    \includegraphics[width=0.49\linewidth]{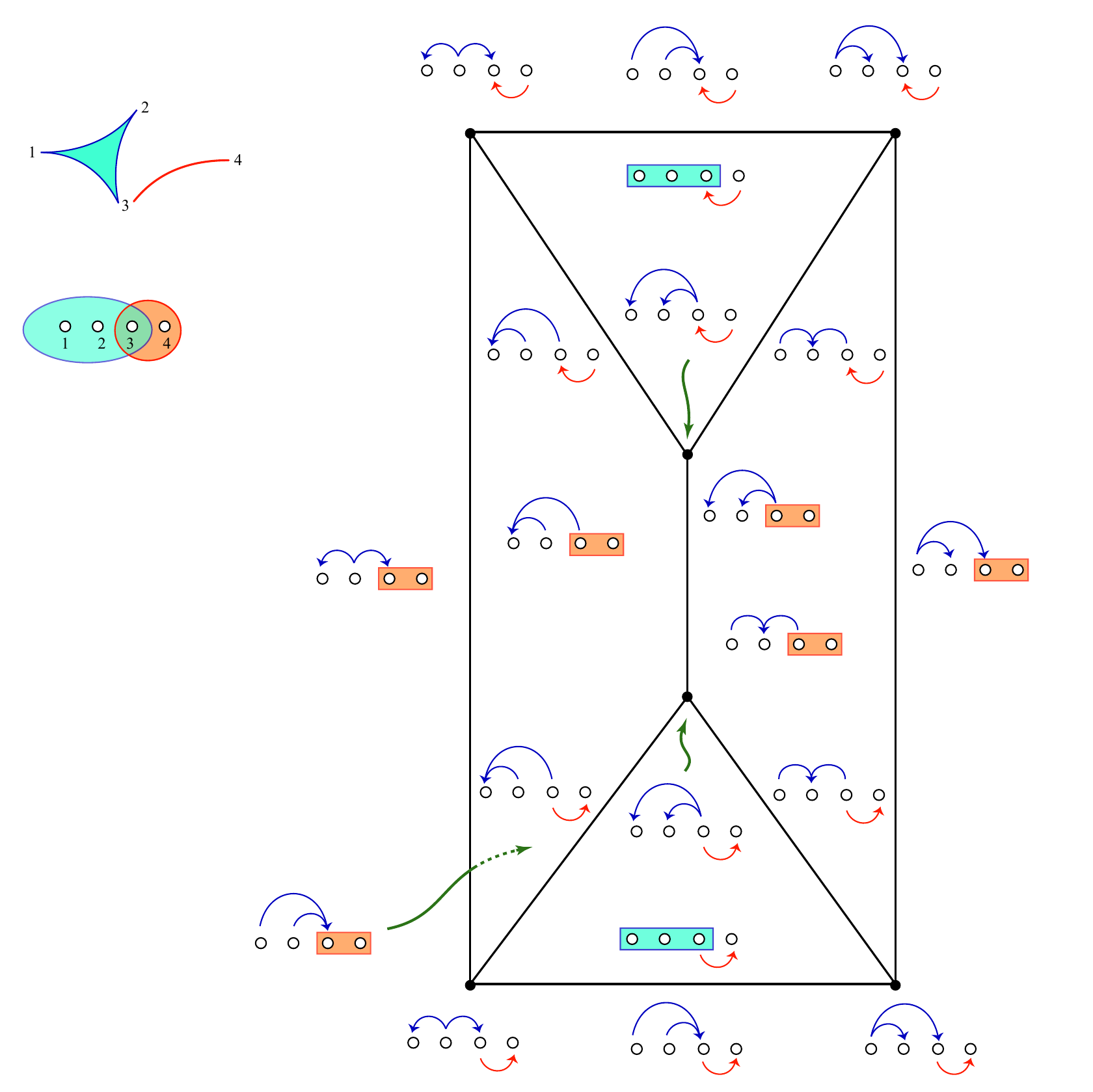}
    \includegraphics[width=0.49\linewidth]{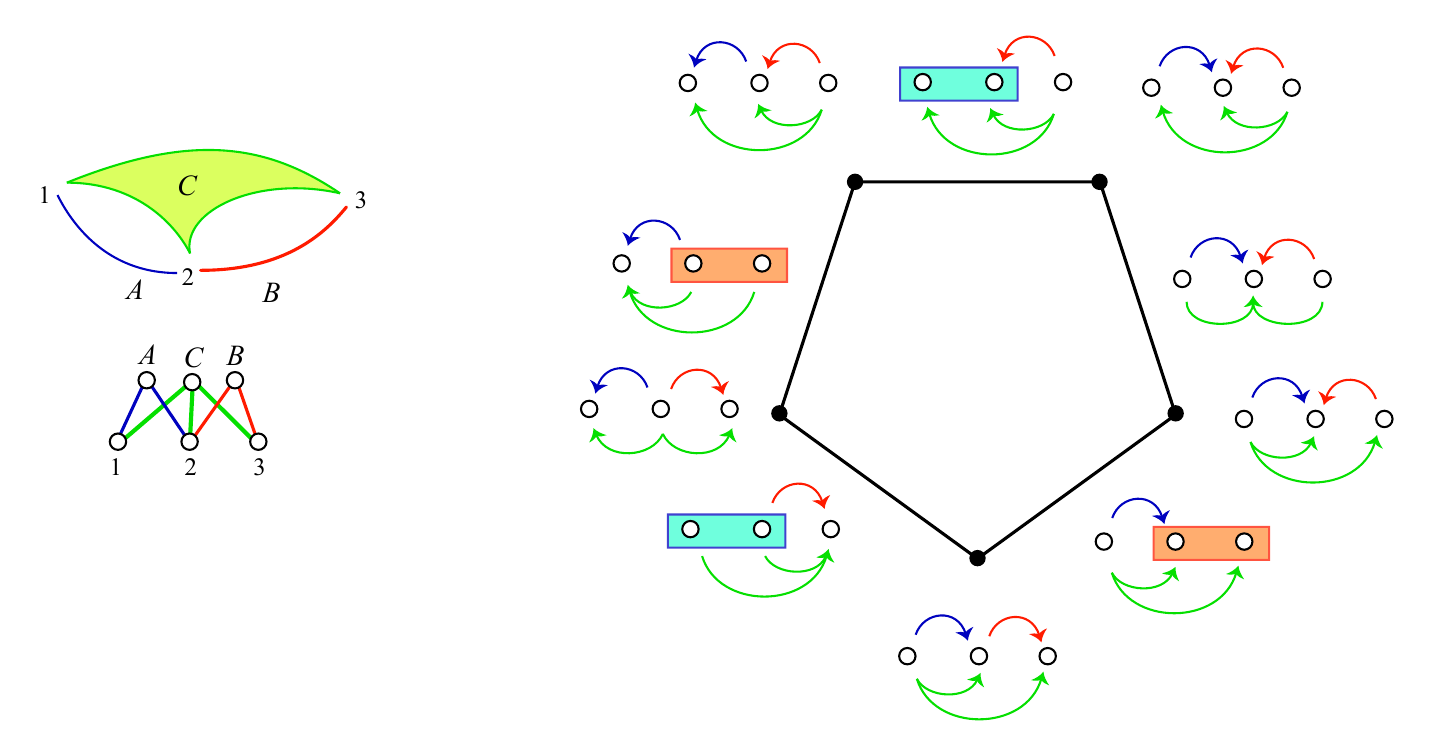}
    \caption{Two examples of hypergraphic polytopes.}
    \label{fig:hypassocjust}
\end{figure}

\begin{defn}  \cite{berg, pilaud_berg}
    The \textit{hypergraphic polytope} $\Delta_G$ for $G$  with $n$ nodes is an $(n-1)$-dimensional polytope; it is the Minkowski sum over unique hyperedges $\displaystyle{\sum_{e\in E}\Delta_e}$ where $\Delta_e$ is the convex hull of standard basis vectors $\mathbf{e}_v \in \mathbb{R}^n$ for $v\in e$. Faces are labeled by the partial acyclic orientations, ordered by refinement. The vertices are the acyclic orientations known as \textit{total orientations}, whose source-sets are all singletons. 
\end{defn}
Examples are in Figures~\ref{fig:hypassocjust} and~\ref{fig:hyp3ord}. The source set of a hyperedge is shown as the set of starting nodes for arrows. When the source set is the entire hyperedge the nodes are highlighted. The hypergraphic polytope is realized as the Minkowski sum in \cite{berg}, and we give a new realization as a cross section of an order polytope in Section~\ref{sec:hyp2}.
\begin{remark}
Our definition of hypergraph allows both repeated hyperedges and hyperedges of cardinality one, which we will refer to as loops. 
The only consistent extension of the definition of acyclic orientations on simple hypergraphs to these multihypergraphs or pseudohypergraphs is that the orientations on two repeated hyperedges (with the same nodes) must be precisely the same. Otherwise the head and tail sets for that same set of nodes would overlap, causing a cycle. Thus the  multi-hypergraphic polytope will be combinatorially the same as the hypergraphic polytope for the underlying simple hypergraph. 
Loops have only one vertex, so it must be the element of the source set, and the target set must be empty. This can be seen as equivalent to saying that loops may not have a direction assigned, since they will be cycles if so.
\end{remark}

\section{Cosmological polytopes}\label{sec:cosmo}

 For a connected pseudograph $G=(V,E)$, the cosmological polytope $\mathcal{C}_G$ has dimension $|V|+|E|-1$ and is defined as the convex hull of its vertices. We follow the definition in \cite{Arkani-Hamed:2}: 
\begin{defn}
 Given a pseudograph $G=(V,E)$, each edge $e\in E$ contributes $|e|+1$ vertices to the cosmological polytope. In particular, an ordinary edge $e=\{u,v\}$ and a loop $\ell=\{w\}$ contribute, respectively,
\begin{align}
    \mathcal{V}_e(\mathcal{C}_G)
    &=
    \{p_e,p_{e,u},p_{e,v}\},
    &
    \mathcal{V}_\ell(\mathcal{C}_G)
    &=
    \{p_\ell,p_{\ell,w}\}.
\end{align}
The vertex set of the cosmological polytope is therefore
$
    \mathcal{V}(\mathcal{C}_G)
    =
    \bigcup_{e\in E}\mathcal{V}_e(\mathcal{C}_G).
$

Each vertex is assigned a point that is found as a vector sum of standard basis elements. For a multigraph $G = (V,E)$ we
choose an arbitrary ordering of the edges $e_1, e_2,\dots, e_m$ and of the  vertices $v_1,v_2,\dots,v_k$. Then we assign each vertex and edge in that order to the standard basis vectors  of $\mathbb{R}^{k+m}$, with subcripts indicating the assignment: $\textbf{e}_{v_1},\dots,\textbf{e}_{v_k},\textbf{e}_{e_1},\dots,\textbf{e}_{e_m}.$
Then for each edge $e=\{u,v\}$, we assign the sums as follows:
  
$p_e \to \textbf{x}_e = \textbf{e}_{u}+\textbf{e}_{v}-\textbf{e}_{e}$

$p_{e,u} \to \textbf{x}_{e,u} = \textbf{e}_{u}-\textbf{e}_{v}+ \textbf{e}_{e}$

$p_{e,v} \to \textbf{x}_{e,v} = -\textbf{e}_{u}+\textbf{e}_{v}+\textbf{e}_{e}$

\noindent Then we define $\mathcal{C}_G$ as the convex hull of the points $\{\textbf{x}_e, \textbf{x}_{e,u}, \textbf{x}_{e,v}\}$ for all edges $e\in E.$  
\end{defn}
\begin{example}
     See the vertex points on the left of Figure~\ref{fig:compare}. Using the coordinate ordering $(x_l,x_e,x_u,x_v),$ we get the points:
$\textbf{x}_l = (-1,0,2,0),~~
\textbf{x}_e = (0,-1,1,1)
,~~ \textbf{x}_{e,u} =  (0,1,1,-1)
,~~ \textbf{x}_{e,v} = (0,1,-1,1)
,\text {and } \textbf{x}_{l,u} = (1,0,0,0,).$  Compare to our new realization, as seen in Example~\ref{newby}.
\end{example}

All the vertex points do indeed turn out to be vertices of the convex hull, as shown in \cite{Arkani-Hamed:2}.\footnote{We are not sure that the current proofs clearly cover the loop edges; if not then the results here settle that question.} One of our goals is to define  $\mathcal{C}_G$ for $G$ a hypergraph, but this classic convex hull realization does not immediately generalize to hyperedges. Our recognition of the combinatorial equivalence of the cosmological polytope to the polar dual of a specific order polytope however will provide geometric realizations for the general hypergraph cosmological polytopes as well as new realizations for the classic case. To prove our main result we will use the following characterizations of the facets and faces of $\mathcal{C}_G$ in terms of subsets of vertices given in \cite{Arkani-Hamed:2}.

\begin{theorem}{\cite{Arkani-Hamed:2}}
Let $G=(V,E)$ be a pseudograph. The facets $F_H$ of $\mathcal{C}_G$ are in bijection with the connected subgraphs $H=(V_H,E_H)$ of $G$, corresponding to the vertex subset $X \subseteq\mathcal{V}(\mathcal{C}_G)$ containing all vertices except: 
\begin{itemize}
    \item $p_e$ for an edge $e\in E_H$;
    \item $p_{e,v}$ for an edge $e\in E$ incident to a vertex $v\in V_H$, with $e\notin E_H$.
\end{itemize}
\end{theorem}
\begin{remark}
The subgraphs $H$ of $G$ need not be induced and include both the single nodes of $G$ and the improper subgraph $H=G$ itself. The facet $F_G$ is called the \textit{scattering facet}.
\end{remark}

\begin{theorem}{\cite{Arkani-Hamed:2}}
Let $G=(V,E)$ be a pseudograph. A subset
$X\subseteq\mathcal{V}(\mathcal{C}_G)$ is the vertex set of a face of
$\mathcal{C}_G$ if and only if it satisfies both of the following
conditions:
\begin{enumerate}
    \item Convexity condition. If for a node $v \in V$ and an edge $e \in E$ incident to $v$, the set $X$ contains both $p_{e}$ and $p_{e,v}$, then, for every edge $f$ incident to $v$, $X$ also contains both $p_f$ and $p_{f,v}$.
    \item Cycle condition. Let $\sigma=(v_1,e_1,v_2,e_2,\ldots, v_k,e_k,v_1)$
    be a cycle of $G$, where $e_i = \{ v_i,v_{i+1}\}$, with indices understood modulo $k$. Then 
    $ \{p_{{e_1,v_{1}}},p_{e_2,v_2}, \ldots p_{e_k,v_1} \}  \subseteq X$ if and only if $\{p_{{e_1,v_{2}}},p_{{e_2,v_{3}}}, \ldots p_{{e_k,v_{1}}} \} \subseteq X$.
\end{enumerate}
\end{theorem}
\begin{remark}
For a loop $\sigma=(v,e,v)$, the two sets appearing in the cycle condition coincide, and the cycle condition therefore imposes no additional constraint.
\end{remark}

\begin{example}\label{uno}
Consider the pseudograph $G$ with vertex set $\{u,v\}$, with loop
$\ell=\{v\}$ and non-loop edge $e=\{u,v\}$, as pictured in Figure~\ref{fig:compare}. The vertex sets corresponding to faces of the cosmological polytope $\mathcal{C}_G$ are
\begin{align*}
\dim 3:\quad&
\{p_\ell,p_{\ell,u},p_e,p_{e,v},p_{e,u}\},
\\[1mm]
\dim 2:\quad&
\{p_\ell,p_e,p_{e,v}\}, \{p_\ell,p_{\ell,u},p_e,p_{e,u}\},
\
\{p_{\ell,u},p_e,p_{e,v}\}, \
\{p_\ell,p_{e,v},p_{e,u}\},
\
\{p_{\ell,u},p_{e,v},p_{e,u}\},
\\[1mm]
\dim 1:\quad&
\{p_{\ell,u},p_{e,u}\},
\
\{p_{\ell,u},p_{e,v}\},
\
\{p_{e,v},p_{e,u}\},
\
\{p_{\ell,u},p_{e}\}, \
\{p_{\ell},p_{e,u}\},
\
\{p_{\ell},p_{e,v}\},
\
\{p_{e},p_{e,v}\},
\
\{p_{e},p_{\ell}\},
\\[1mm]
\dim 0:\quad&
\{p_{\ell,u}\},
\
\{p_{\ell}\},
\
\{p_{e,u}\},
\
\{p_{e,v}\},
\
\{p_{e}\},
\\[1mm]
\dim(-1):\quad&
\varnothing.
\end{align*}
As $G$ contains no non-loop cycles the only condition is the convexity condition at vertex $u$ which forbids subsets which contain exactly one of $\{ p_\ell, p_{\ell,u}\}$ or $\{ p_e, p_{e,u}\}$.
\label{ex:one}
\end{example}

\section{Order polytopes}\label{sec:order}
For a finite poset $P$, we follow Galashin \cite{galashin} and consider the order polytope $\mathcal{O}(P)$ of dimension $|P|-2.$ (See Remark~\ref{stan} for comparison with Stanley's order polytope.) The face lattice of the order polytope is described by a subset of partitions of $P$ which we now describe.

A partition $\pi=\pi_1| \ldots |\pi_k$ of $P$ is a collection of non-empty, pairwise disjoint subsets $\pi_i\subseteq P$, called \emph{parts}, whose union is $P$. Given $p\in P$, we denote by $\pi_p$ the unique part of $\pi$ containing $p$. When $P$ is equipped with a partial order, it is useful to introduce the following terminology.
\begin{defn} 
A partition $\pi$ of a poset $P$ is called \emph{connected} if, for every part $\pi_i\in\pi$, the subgraph of the Hasse diagram of $P$ induced by $\pi_i$ is connected. 
\end{defn}

\begin{defn} Given a partition $\pi$ of a poset $P$, we define the directed graph $P/\pi$ whose nodes are the parts of $\pi$ and directed edge \( \pi_i\rightarrow\pi_j \) whenever there exist $p\in\pi_i$ and $q\in\pi_j$ such that $q$ covers $p$ in $P$.
\end{defn} 

\begin{defn} 
The partition $\pi$ of $P$ is said to be \emph{acyclic} if $P/\pi$ contains no directed cycles. We denote by $\Pi(P)$ the set of connected, acyclic partitions of $P$. 
\end{defn} 

\begin{remark} If $\pi\in\Pi(P)$, then every part of $\pi$ is order-convex, that is, if $p,r\in\pi_i$ and \( p<q<r \) for some $q\in P$, then $q\in\pi_i$.  
\end{remark} 

\begin{defn} For a finite poset $P$, the order polytope $\mathcal{O}(P)$ is a convex polytope of dimension $|P|-2$ whose face lattice is isomorphic to $\Pi(P)$ ordered by reverse refinement. We review geometric realizations of the order polytope in Section~\ref{root}.
\end{defn}

Under this correspondence, the facets of $\mathcal{O}(P)$ are labeled by partitions with a unique non-singleton part of the form \( \{p,q\}\) where $q$ covers $p$ in $P$, an edge of the Hasse diagram of $P$. The vertices of $\mathcal{O}(P)$ are labeled by two part partitions $\pi \in \Pi(P)$. Note that for a poset with a unique maximal element, the vertices will have one part of their two part partition containing that max, and by convexity that part will be an upper ideal of the poset.  We will primarily work with the polar dual polytope $\mathcal{O}(P)^{\Delta}$, whose face lattice is dual to that of $\mathcal{O}(P)$. In particular,
vertices and facets are exchanged. Thus in Figure~\ref{fig:compare} the facets of $\mathcal{O}(P'(G))^{\Delta}$ are labeled by two part partitions $\pi \in \Pi(P)$.

\begin{example}\label{dos}
Consider again the graph of example $\ref{ex:one}$. As pictured in Figure~\ref{fig:compare}, the connected, acyclic partitions labeling faces of  $\mathcal{O}(P'(G))^{\Delta}$ are:
\begin{align*}
\dim 3:\quad& \ell u e v \star,
\\[1mm]
\dim 2:\quad&
u|\ell e v \star, 
\
v | \ell u e \star,
\
\ell u| e v \star, 
\
u e v| \ell \star,
\
\ell u e v| \star,
\\[1mm]
\dim 1:\quad&
\ell u e|v|\star ,
\
\ell u|e v|\star ,
\
\ell|u e v|\star,
\
\ell u|v|e  \star, \
u e|v|\ell  \star,
\
u|e v|\ell \star,
\
\ell|u|e v \star,
\
u|e|v|\ell e  \star ,
\\[1mm]
\dim 0:\quad&
\ell u|e|v|\star,
\
u|e|v| \ell \star ,
\
\ell|u e|v|\star,
\
\ell|u|e v|\star,
\
\ell|u|v| e \star,
\\[1mm]
\dim(-1): \quad & \ell|u|e|v| \! \star.
\end{align*}
\end{example}

\begin{remark}\label{stan}
    Order polytopes $\mathcal{O}_{St}(P)$  (with subscript added for clarity) were first defined by Stanley in \cite{stanley} as the set of monotonic functions from a poset $P$ to the unit interval. That is, the output of the function for each element $i\in P$ is seen as a variable $0\le x_i \le 1$, and the list of outputs from a function  as coordinates of a point $\mathbf{x} = (x_1,x_2,\dots,x_n).$ These points make a 0/1 polytope of dimension $n,$ inside the unit hypercube $[0,1]^n$.  Equivalently, we can adjoin a new unique min and max to $P$ often called $\{\hat{0},\hat{1}\}$ respectively. The new \textit{bounded} poset is called $\hat{P}$ and the polytope $\mathcal{O}_{St}(P)$ is affinely (and thus combinatorially) equivalent to a new version called $\hat{\mathcal{O}}_{St}(\hat{P})$, still $n$-dimensional but now living in $n+2$ dimensions. $\hat{\mathcal{O}}_{St}(\hat{P})$ is defined by setting the coordinates corresponding to the new min and max equal to 0 and 1 respectively, while keeping all the inequalities $x_i\le x_j$ for $i\precdot j$ a covering relation in $\hat{P}.$

    The polytope we call $\mathcal{O}(P)$, given by the connected acyclic partitions of $P$, is described first by Galashin in \cite{galashin}, where it is shown to be combinatorially equivalent to Stanley's $\hat{\mathcal{O}}_{St}(P)$ when $P$ is bounded, with unique min and max. Thus Galashin's order polytope is closely related to Stanley's: if $P$ is any poset, then forming $\hat{P}$ will mean that Galashin's $\mathcal{O}(\hat{P})$ is equivalent to Stanley's $\hat{\mathcal{O}}_{St}(\hat{P})$ which is equivalent to Stanley's $\mathcal{O}_{St}(P).$ On top of that, for any connected $P$, the polar dual of Galashin's $\mathcal{O}(P)$ is a codimension-2 face of the polar dual of Galashin's $\mathcal{O}(\hat{P})$. That's because faces are indexed by connected partitions, and $P$ itself is part of the partition of $\hat{P}$ into three parts: the max, the min and $P.$ 
    But for $\hat{P}$, Galashin's and Stanley's versions are combinatorially equivalent, so the polar dual of Galashin's $\mathcal{O}(P)$ is a face of Stanley's version as well.

\end{remark}

\section{Combinatorial equivalence}\label{sec:equivo}
In this section we show that the cosmological polytope $\mathcal{C}_G$ is combinatorially equivalent to the polar dual of the order polytope $\mathcal{O}(P'(G))^{\Delta}$ for a poset $P'(G)$ constructed from the incidence structure of $G$.

\begin{figure}
    \centering
    \includegraphics[width=\linewidth]{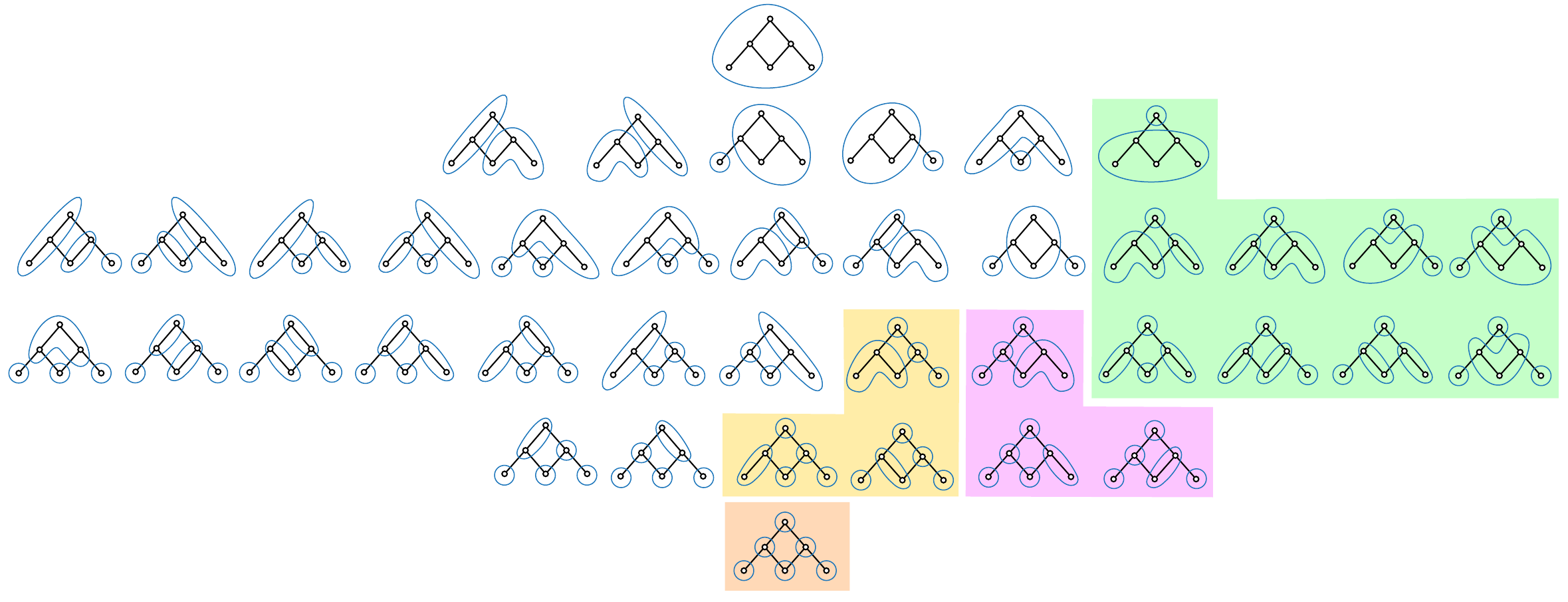}
    \caption{Faces of $\mathcal{C}_G \cong \mathcal{O}(P'(G))^{\Delta},$  for $G$ the path graph with two edges. Highlighted on the right in four colors is the scattering facet $F_G$, a tetrahedron partitioned into the face posets of hypergraphic polytopes $\Delta_G$ and its edge-deleted subgraph hypergraphic polytopes. The largest of these highlighted face posets is pictured as the hypergraphic polytope of the path graph in Figure~\ref{fig:hyppathjust}.}
    \label{fig:cosmo-order}
\end{figure}

\begin{defn} 
Let $G=(V,E)$ be a pseudograph. The \emph{incidence poset} of $G$ is the rank-2 graded poset $P(G)=V \sqcup E$ with cover relations $v\lessdot e$ whenever $v$ is incident to $e$ in $G$. The \emph{augmented incidence poset} $P'(G)=V\sqcup E\sqcup\{\star\}$ is obtained from $P(G)$ by adjoining a maximal element $\star$, with cover relations $e\lessdot\star$ for every $e\in E$.
\end{defn}
\begin{example}
The poset $P'(G)$ constructed from the graph with one edge and one loop is displayed in Figure~\ref{fig:compare}. 
The poset $P'(G)$ constructed from the path graph with two edges is displayed in Figure~\ref{fig:cosmo-order}. 
\end{example}

\begin{figure}
    \centering
    \includegraphics[width=\linewidth]{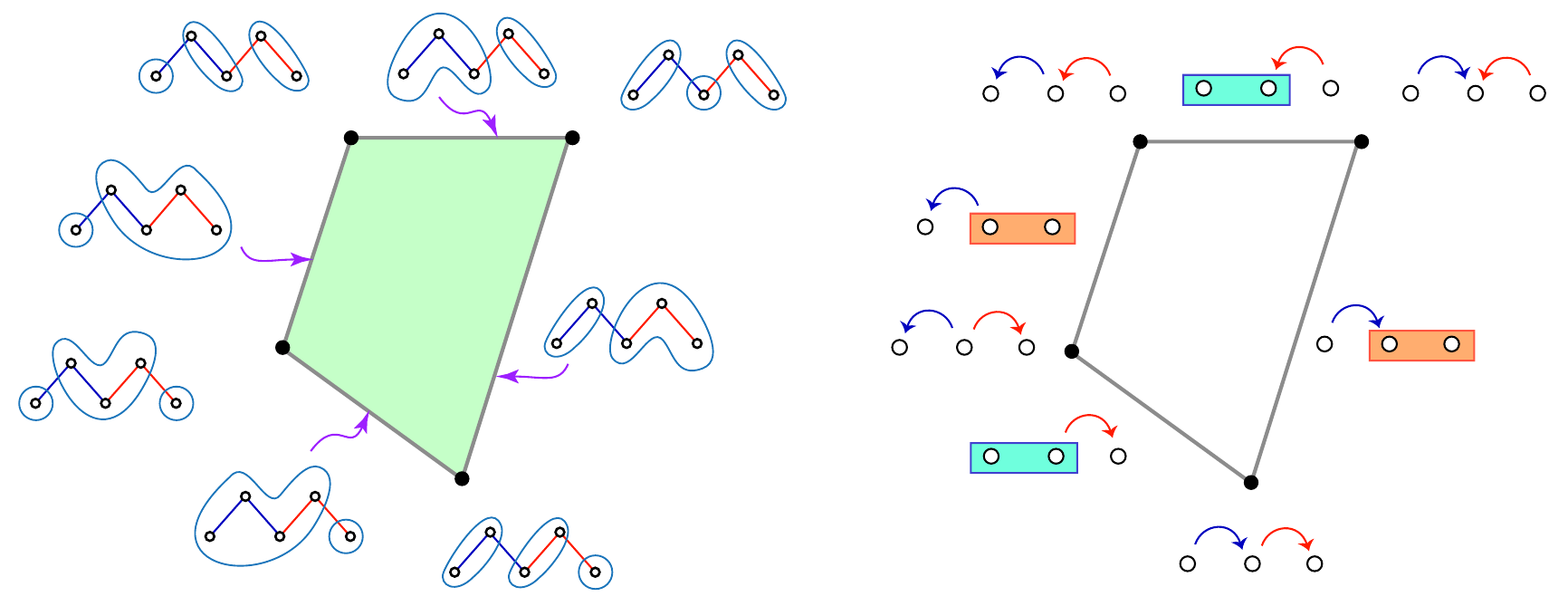}
    \caption{Left to right, the bijection described in the proof of Theorem~\ref{conj}. On the right: the hypergraphic polytope $\Delta_G$ for the path graph with two edges. On the left: the same face poset seen as acyclic partitions.  This quadrilateral occurs as a cross section of the tetrahedral scattering facet with faces highlighted in Figure~\ref{fig:cosmo-order}: the acyclic partitions are seen there with the singleton maximal element added on.}
    \label{fig:hyppathjust}
\end{figure}

\begin{theorem}\label{big}
    Let $G=(V,E)$ be a connected pseudograph. The cosmological polytope $\mathcal{C}_G$ is combinatorially equivalent to the polar dual order polytope $\mathcal{O}(P'(G))^{\Delta}.$
\end{theorem}

\begin{proof}
    The dimension of $\mathcal{C}_G$ is $|V|+|E|-1,$ which is equal to $|P'(G)|-2$, the dimension of $\mathcal{O}(P'(G))^{\Delta}.$ For vertex-facet incidence to match between the polytopes we need bijections, between vertices and facets, and for vertex incidence to be preserved under that correspondence. We define a bijection on vertices between the two polytopes by taking an edge of the Hasse diagram of $P'(G)$ to a vertex of $\mathcal{C}_G.$ Recall that each element in the second rank of $P'(G)$ corresponds to an edge $e\in E$.  An edge of the Hasse diagram between the maximal element $\star$ and $e$ is mapped to $p_e$. An edge of the Hasse diagram betwen $e$ and $v$ for $v\in V$ is mapped to $p_{e,v}$.

    For facets, the bijection is to take a facet of the order  polytope $\mathcal{O}(P'(G))^{\Delta}$, which is a 2-part partition, to the connected subgraph $H=(V_H,E_H)$ of $G$ which corresponds via the edge and node elements to the part of that partition \textbf{not} containing the maximal element $\star.$

    Since both our structures are known to be polytopes, now we need only check that our bijection and its inverse preserve vertex-facet incidence. We claim that if an edge of the Hasse diagram is a subset of a part of a 2-part partition (labelling a facet of $\mathcal{O}(P'(G))^{\Delta}$), then its image $p_e$ or $p_{e,v}$ will be a member vertex of the facet of $\mathcal{C}_G$ labeled by the part of that partition not containing the maximal element $\star.$ This is clear since if the edge of the Hasse diagram is from $\star$ to the edge-element $e$ then the edge $e$ of the graph is not in the subgraph $H$ corresponding to the part not containing $\star$. Also, if the edge of the Hasse diagram is from element $e$ to element $v$, either the node $v$ is not in the subgraph $H$ (if the element $v$ is in the part containing $\star$) or the node $v$ is in $H$ but the edge $e$ is also in $H.$   The preservation of vertex facet incidence is also easily checked in the inverse direction of the bijection.
 \end{proof}

An example of the bijection is seen in Figure~\ref{fig:compare}. 


\subsection{Generalized definition}
Since the cosmological polytopes are seen as special versions of order polytopes, we can immediately generalize the definition of cosmological polytope of a pseudograph: to hypergraphs and to multi and pseudo hypergraphs, that is, to multisets of elements from any power set of vertices. Combinatorially, the cosmological polytope of any pseudo-hypergraph $G$ is the dual of the order polytope of the poset: $ C_G := \mathcal{O}(P'(G))^{\Delta}.$ Via the realizations of the order polytope there are new geometric realizations available for the cosmological polytopes, although they may not be directly interesting to physicists.

Given a connected pseudo-hypergraph $G$, the scattering facet of the cosmological polytope is the facet labeled by $G$ itself. Seeing the cosmological polytope as the polar dual of the order polytope of the augmented incidence poset $P'(G),$ the scattering facet is labeled by the partition of $P'(G)$ one part of which is the unique maximal element and the other part of which is the incidence poset $P(G)$, with hyperedges as maximal elements covering the nodes they contain. Thus we have:

\begin{corollary}\label{scatter}
    The scattering facet of a  cosmological polytope for a connected pseudo-hypergraph $G$ is combinatorially equivalent to the polar dual of the Galashin order polytope of the (unaugmented) incidence poset of that pseudo-hypergraph: $F_G \cong \mathcal{O}(P(G))^{\Delta}$.  
\end{corollary}
 Even further, we define the cosmological polytope of any poset $P$ as simply $C_P = \mathcal{O}(P\sqcup  \star)^{\Delta},$ the polar dual of the order polytope of  $P$ augmented with the new unique maximal element.

 \section
   {New realizations for cosmological polytopes.}\label{root}

A \textit{root polytope} is a convex hull of some standard basis vectors, some negatives of standard basis vectors, and some differences of standard basis vectors. A connected \textit{starred quiver} $Q$ is a connected directed graph with sources and sinks (together called the starred set) and no edge from a source to a sink.  In \cite{williams} there is  defined a root polytope called Root($Q$) which is the convex hull of points, based on standard basis vectors indexed by the non-starred elements of $Q.$ Each edge of $Q$ gives rise to a point: edges $e$ from a source to a node $u$ give rise to the point $\textbf{e}_{u}$, edges from $u$ to a sink the point $-\textbf{e}_{u}$, and other (non-starred) edges directed from $u$ to $v$ the difference $\textbf{e}_{v}-\textbf{e}_{u}.$ 

Theorem D of \cite{williams} compares the polar dual of Stanley's order polytope with the root polytope of a quiver based on the Hasse diagram of $\hat{P}$, where graded $P$ is augmented with both a min and max.  We restate just the part that applies to our posets $P(G)$:
\begin{theorem}\cite{williams}
Let $P$ be a finite graded poset and let $Q:=Q_{\hat{P}}$ be the associated starred quiver: an upwards arrow for each edge of the Hasse diagram of bounded $\hat{P} = P\sqcup\{\hat{0},\hat{1}\}.$   Then the face fan $\mathcal{F}_Q$ 
of the root polytope $\text{\emph{Root}}(Q)$ of the starred quiver $Q$ coincides with the (inner) normal fan $\mathcal{N}(O_{St}(P))$. 
\end{theorem}
Since the vertex-and-edge incidence poset $P(G)$ for a hypergraph $G$ is graded, Theorem D implies that the inner normal fan of the root polytope is equal to the face fan of the root polytope of the starred quiver of the poset $P(G)$ with both an adjoined max and adjoined min. That is, in the language of this paper:
\begin{corollary}
    For a pseudo-hypergraph $G$, the root polytope $\text{\emph{Root}}(Q_{\hat{P}(G)})$ has a facet that is combinatorially equivalent to the cosmological polytope $\mathcal{C}_G.$ Therefore, the cosmological polytope is itself a root polytope.
\end{corollary}
\begin{proof}
 Corollary 4.19 of \cite{williams} points out that the root polytope  in this case will be combinatorially equivalent to the polar dual of the order polytope of Stanley.   
    The facet of that polar dual we are interested in is the one which corresponds to the partition into two parts : $\hat{P}(G)= P'(G)\sqcup \{\hat{0}\}$, our augmentation and the singleton unique min. It is thus equivalent to the cosmological polytope by Theorem~\ref{big}. Of course, any face of a root polytope is by definition a root polytope itself. \end{proof}
    
    Furthermore, being a combinatorial copy of $\mathcal{C}_G$,  this facet  has its vertices corresponding to the edges of the Hasse diagram of $P'(G)$. A set of vertices of the root polytope that has that same number of vertices, and enough information to determine the Hasse diagram of $P'(G)$, is the set found by excluding  from the root polytope $\text{Root}(Q_{\hat{P}(G)})$ the vertices (standard basis vectors) that involve the unique min.
    \begin{corollary}\label{newg}
        The vertices of $\text{\emph{Root}}(Q_{\hat{P}(G)})$ which correspond to the edges of the Hasse diagram of $P'(G)$ are the vertices of a facet combinatorially equivalent to $\mathcal{C}_G$.
    \end{corollary}
    \begin{proof}
        In the proof of Theorem 4.15 in \cite{williams} the facets of $\text{Root}(Q_{\hat{P}(G)})$ are shown to correspond to certain pairs of connected subquivers, one containing the max and the other the min. These match the connected acyclic partition of two parts labeling the facets of the polar dual of Stanley's order polytope, and the facet which is the polar dual of Galashin's order polytope is the one we desire: the partition into $P'(G) \sqcup {\hat{0}}$. The facet equality using coordinate variables indexed by the nodes of $P(G)$,  according to \cite{williams}, will be $\sum_e x_e +2\sum_u x_u \ge -1.$ This will be satisfied (as equality) by the points corresponding to edges in $P'(G)$, specifically $-\textbf{e}_{e}$ and $\textbf{e}_{e}-\textbf{e}_{u},$ and strict inequality for $\textbf{e}_{u}.$ 
    \end{proof}
    Therefore, in order to realize $\mathcal{C}_G$, we list the edges of $G$ and the nodes of $G$  in arbitrary order as indices of standard basis vectors. Then we assign each covering relation of $P'(G)$ to a combination of those standard basis vectors. The new vertices of the cosmological polytope we now call $q_e$ and $q_{e,u},$ since each hyperedge and each of its member nodes contributes a point. 
    
    $q_e \to \textbf{y}_e = -\textbf{e}_{e}$,

$q_{e,u} \to \textbf{y}_{e,u} =\textbf{e}_{e} -\textbf{e}_{u}$,

\noindent and the cosmological polytope is found as the convex hull of all these $\textbf{y}_{e}$ and $\textbf{y}_{e,u}.$

\begin{example}\label{newby}
    See the vertex points on the left of Figure~\ref{fig:compare}. Using the coordinate ordering  $(y_l,y_e,y_u,y_v),$ we get the points:
$\textbf{y}_l = (-1,0,0,0),~~
\textbf{y}_e = (0,-1,0,0)
,~~ \textbf{y}_{e,u} = (0,1,-1,0) 
,~~ \textbf{y}_{e,v} = (0,1,0,-1)
,\text {and } \textbf{y}_{l,u} = (1,0,-1,0).$
\end{example}

\subsection{Order polytope realizations}

 Galashin describes a realization for $\mathcal{O}(P)$  for real variables $x_i$ indexed by the $n$ elements $i\in P$. This implies a realization for the cosmological polytope $\mathcal{C}_G \cong \mathcal{O}(P'(G))^{\Delta}$, for $G$ any (pseudo)hypergraph:    
 \begin{enumerate}
     \item facet inequalities: $x_i\le x_j$ for each covering relation (edge of the Hasse diagram) $i\precdot j$ in $P'(G)$,
     \item equalities: $\displaystyle{\sum_{i\in P'(G)}x_i = 0}$ and  $\displaystyle{\sum_{i \precdot j}(x_j-x_i) = 1,}$
     \item and then find the polar dual.
 \end{enumerate}
There is another realization for $\mathcal{O}(P)$   in \cite[Sec. 8.2]{mppfull} for  variables $x_{e}$ indexed by the $|E|$ edges of the Hasse diagram of $P$ (covering relations of $P$), as a cross section of a simplex. That allows us to find the cosmological polytope by a polar dual to that cross section, with coordinates indexed by incidences of vertices in hyperedges:
 \begin{enumerate}
     \item inequalities: $x_{e} \ge 0$ for each edge $e$ of the Hasse diagram of $P'(G)$,
     \item equality: $\displaystyle{\sum_e x_{e} = 1,}$
     \item equalities for independent cycles $c$ in the underlying graph of the Hasse diagram of $P'(G)$: Choosing an arbitrary direction to traverse the cycle $c$, the orientation (covering relation) of the edge $e$ in $c$ will either agree with or be against the direction of traversal of $c$, and the equality is $\displaystyle{\sum_{e \text{ with }c }x_e = \sum_{e \text{ against } c}x_e },$  
     \item and then find the polar dual.
 \end{enumerate}


Notice that the last two realizations described here end with finding the polar dual, which means the specific inequalities and vertices of the polytope will be determined by the algorithm used. It would be interesting to know if one of these actually recovers the root polytope realization of Corollary~\ref{newg}. Any of the realizations here can be restricted to find a realization of the scattering facet $F_G$, which is combinatorially equivalent to $\mathcal{O}(P(G))^{\Delta}$ by Corollary~\ref{scatter}.

\section{Hypergraphic polytopes as cross sections of dual order polytopes.}\label{sec:hyp2}

We demonstrate that given a connected hypergraph $G$, the hypergraphic polytope $\Delta_G$ is found as a cross section of the polar dual of the order polytope of the poset $P(G)$. 
   Our explicit map of faces is new, but the geometric realization can be seen as a specialization of results of \cite[Sec. 14]{postnikov}, which uses root polytopes of bipartite graphs.  

To begin, we prove a lemma: given a connected hypergraph $G$, the polar dual of the order polytope of the poset $P(G)$  contains a model of the hypergraphic polytope $\Delta_G$. This can be described as a new model for the faces of the hypergraphic polytope $\Delta_G$. The faces of the hypergraphic polytope can be labeled by all the acyclic partitions of the poset $P(G)$ which have no parts that are a single maximal element (no singleton parts whose element is a hyperedge.) Precisely, we can say the following, where we ignore multiedges and loops (since they do not contribute to the hypergraphic polytope.)

\begin{lemma}\label{hyper}
  For a simple connected hypergraph $G$, the face poset of the hypergraphic polytope $\Delta_G$ is found as a facial subcomplex of the polar dual of the order polytope, $\mathcal{O}(P(G))^\Delta$ for the (unaugmented) incidence poset $P(G)$. 
\end{lemma}

\begin{proof}
    We describe a bijection from our selected acyclic partitions of the poset of $G$ to the acyclic orientations of $G.$ Given an acyclic partition of the poset $P(G)$ with no singleton second rank parts, we describe the source-sets (head sets) of the hyperedges of $G$ as follows: for each hyperedge of $G$ the source-set is the set of nodes in that hyperedge which is a subset of the part of the partition that contains that hyperedge itself (as a maximal node of the rank 2 poset). 

    The restriction to acyclic partitions that have no part which is a singleton consisting of only a maximal node ensures that the head set of each hyperedge is non-empty. The acyclic requirement for the partition of the poset translates directly to the acyclic requirement for the orientation of the hypergraph.

    Covering refinement of the partition of the poset (avoiding any singletons of second rank) implies covering refinement of the acyclic orientation and vice versa.  
Under our bijection the acyclic partitions of $P(G)$ into two parts will map to facets of $\Delta_G$.  Thus our bijection exhibits the face poset of $\Delta_G$ as a facial subcomplex of the $\mathcal{O}(P(G))^\Delta$.
 \end{proof}

\begin{remark}
    For a pseudo-hypergraph $G$ with at least one unique node for each edge, all the two-part partitions of the poset will be the kind we want in our complex: no singleton hyperedge (second rank) parts. Notice that in this case the two-part partitions label both the all the facets of the scattering facet $F_G$ and the facets of the hypergraphic polytope $\Delta_G$. Thus in this case the complex includes all the facets of the polar dual of the order polytope, i.e. all the vertices of the order polytope, so it constitutes a \textit{spanning} complex of $\mathcal{O}(P(G))$.
\end{remark} 
 
 An example is shown in Figure~\ref{fig:hyp3ord}. Another example is shown in  Figures~\ref{fig:cosmo-order} and~\ref{fig:hyppathjust}.  Now we can demonstrate the geometric meaning:
\begin{theorem}\label{conj}
     For any pseudohypergraph $G$, the hypergraphic polytope $\Delta_G$ can be realized as a cross section of $\mathcal{O}(P(G))^{\Delta}$. 
\end{theorem}
\begin{proof}
    For $G=(V,E)$ we have  by Lemma ~\ref{hyper} 
    a labeling of the hypergraphic polytope with acyclic partitions of $P(G)$, with no singleton part of second rank, ordered by refinement. We use the root polytope geometric realization of the scattering facet $F_G = \mathcal{O}(P(G))^{\Delta}$, using $|E|+|V|$ real variables $x_1,x_2,\dots,x_n$. Thus $F_G$ is the convex hull of the points $\mathbf{x}_{e,u} = \mathbf{e}_e - \mathbf{e}_u$, one for each covering relation in $P(G).$ These points correspond to the vertices of $F_G,$ which are the partitions with one part one edge of the Hasse diagram  and the rest of the parts singletons. 

    Now we take our realization of $F_G = \mathcal{O}(P(G))^{\Delta}$  and we intersect it with the subspace formed by the intersection of hyperplanes $\displaystyle{x_e = \frac{1}{|E|}}$, one for each $e\in E.$ That is, we claim that the hypergraphic polytope is given by the subset of the scattering facet with coordinates corresponding to edges $e \in E$ all set equal to that reciprocal value:
    $$\Delta_G = F_G\cap\mathcal{E}, \text{~~~~where~~~} \mathcal{E} = \{\mathbf{x}\in \mathbb{R}^n \mid \forall e\in E, x_e  = \frac{1}{|E|} \}.$$

    To prove the claim, we use the fact that the face poset structure of the hypographic polytope is determined by refinement of our acyclic partitions of $P(G).$ Then we just need to show that the relative interiors of faces of $F_G$ labeled by acyclic partitions with no second rank singletons do indeed intersect  $\mathcal{E},$ while any face labeled by an acyclic partitions with second rank singletons fails to intersect $\mathcal{E}.$ Then the intersecting faces will inherit the correct membership structure from $F_G$. 

    In the case that a face  of $F_G$ is labeled by an acyclic partition with no second rank singletons, then it will include at least one vertex $\mathbf{y}_e$ of $F_G$  for each second rank element $e$ of the Hasse diagram of $P(G).$  We can pick that vertex just by keeping an edge of the Hasse diagram incident on $e$ that is in a part of the acyclic partition and refining the rest of the parts into singletons. Those vertices will each have $x_e = 1$, and thus the centroid of those points will lie in  $\mathcal{E},$ ensuring that the relative interior of our face intersects $\mathcal{E}.$ 

    For the case that a face  of $F_G$ is labeled by an acyclic partition with a second rank singleton $e$, we note that every vertex of that face will have that coordinate $x_e = 0.$ Thus it will fail to intersect $\mathcal{E},$ and we have shown our claim. 
\end{proof}

\begin{example}
    For the hypergraph  $G$ in Figure~\ref{fig:hyp3ord}, consider the full scattering scattering facet $F_G$ shown in Figure~\ref{fig:33slice}. The hypergraphic polytope is found as the intersection of the convex hull of the six points: $(1,0,-1,0,0,0),(1,0,0,-1,0,0),(1,0,0,0,-1,0)$ for the first hyperedge, and \\$(0,1,0,-1,0,0),(0,1,0,0,-1,0),(0,1,0,0,0,-1)$ for the second, with the space $\mathcal{E}$ given by two equations: $x_1 = 1/2$ and $x_2=1/2$ where $x_1$ and $x_2$ are the first and second coordinates, corresponding to the first and second hyperedges.
\end{example}

\subsection{Partition of scattering facet face poset into hypergraphic polytopes}
  
An \textit{edge-deleted sub-hypergraph} $H\subseteq G$ is formed by keeping all the vertices of $G$ but deleting some hyperedges.

\begin{landscape}
 \begin{figure}
    \centering
    \includegraphics[width=\linewidth]{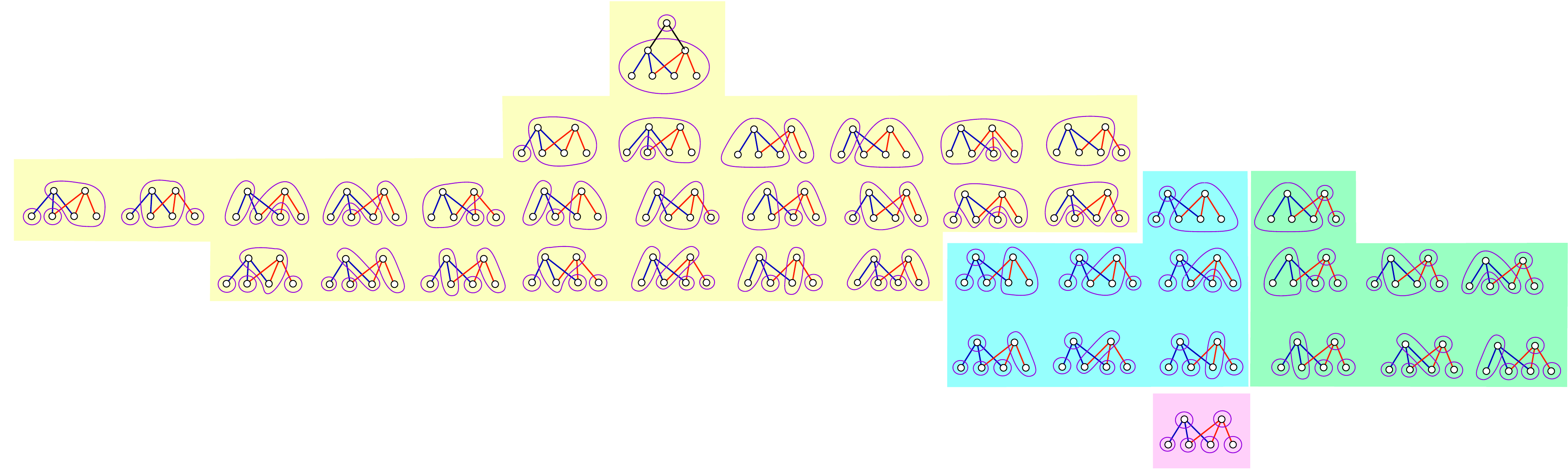}
    \caption{The faces of the scattering facet $F_G$ of $\mathcal{C}_G$ for the hypergraph $G$ from Figure~\ref{fig:hyp3ord}; equivalently, the faces of the polar dual of the order polytope of the incidence poset $P(G)$ of that hypergraph. (The very uppermost Hasse diagram shows the augmented  $P'(G)$  but the rest show only the poset $P(G)$ itself.) Highlighted subposets of faces show the face posets of hypergraphic polytopes for edge-deleted sub-hypergraphs, with $\Delta_G$ at the top, as pictured in Figure~\ref{fig:hyp3ord}. }
    \label{fig:33slice}
\end{figure}
\end{landscape}

\begin{corollary}
    Given a multi-hypergraph $G$, the polar dual of the order polytope of the poset of $G$  admits a partition of its faces into subposets which are themselves the face posets of the hypergraphic polytopes $\Delta_H$ of the edge-deleted sub-hypergraphs $H\subset G.$ 
\end{corollary}
\begin{proof}
    To a sub-hypergraph formed by deleting some hyperedges from $G$ we associate the partitions of the poset of $G$ with those hyperedges as singleton parts. Then the bijection in the proof of Lemma~\ref{hyper} is extended, taking such a partition to an acyclic orientation of the sub-hypergraph. Figures~\ref{fig:cosmo-order} and~\ref{fig:33slice} show the face posets of the hypergraphic polytopes of edge-deleted sub-hypergraphs as highlighted portions 
\end{proof}

Therefore we can describe the $f$-vector of the the scattering facet for $G$ as a sum of $f$-vectors for the hypergraphic polytopes of the edge-deleted (keeping all nodes) sub-hypergraphs $H$ of $G$.  We do so by means of the $f$-polynomial of a polytope $\mathcal{P}$, which is defined:
$$f_{\mathcal{P}}(t) = \sum_{i=0}^{dim(\mathcal{P})} t^{f_i}, \text{ where } f_i \text{ is the number of faces of dimension } i $$ We can also define $\hat{f}(t) = t^{-1}+f(t),$ in order to accommodate the empty face.  

\begin{corollary}\label{effy}
    The $f$-polynomial of the scattering facet ${F}_G$ for $G=(V,E)$ is just the shifted sum of the $f$-polynomials of the hypergraphic polytopes $\Delta_H$ of the edge-deleted sub-hypergraphs $H = (V,E_H)$, including the whole hypergraph $G$. 
    $$\hat{f}_{{F}_G}(t) = \sum_{H \subseteq G} t^{|E_H|-1}f_{\Delta_H}(t), \text{ \emph{ summed over } } H=(V,E_H).$$
\end{corollary}



\bibliographystyle{alphaurl}
\bibliography{hyper}

\section*{AI Disclosure}

AI tools were largely set aside in the production of this work, partly due to the fact that they were tried and found wanting for the purpose. That criticism might apply mostly to the older versions that we queried during exploratory phases: as far as we can recall none of the answers turned out to be directly helpful, although they might have inspired us to think harder via reverse psychology. The actual computations, patterns we found, theorems we proved, reasoning in the proofs, the writing of the text and creation of the figures---were all done ``by hand,'' using standard tools like Overleaf, Polymake, and Adobe Illustrator.

There was one instance of helpful AI use. Proofreading the final draft with AI (GPT-6 Astra) turned up several small errors: one faulty statement we deleted, several that only needed an added qualifier and one suggested attribution. For the latter the AI helped confirm an authorial suspicion that our cross-sectional realization of hypergraphic polytopes was a special version of a construction in Postnikov's previously published paper \cite{postnikov}. We added an attribution to the prior work, at the beginning of Section~\ref{sec:hyp2}.

\end{document}